\documentclass[12pt, reqno]{amsart}

\usepackage[T1]{fontenc}
\usepackage{times}
\usepackage{color,graphicx}
\usepackage{subfig,array}
\usepackage{amsthm,amssymb,amsmath}
\usepackage{hyperref}\usepackage{matlab-prettifier}
\usepackage{pgf,tikz}
\usepackage{mathrsfs}
\usetikzlibrary{arrows}

\calclayout

\newcommand\ket[1]{\ensuremath{|#1\rangle}}
\newcommand\bra[1]{\ensuremath{\langle#1|}}
\newcommand{\braket}[2]{\langle #1 | #2 \rangle}
\newcommand{\ketbra}[2]{| #1 \rangle\langle #2 |}

\def\H{\mathcal{H}}

\newcommand{\bes}{\begin{eqnarray*}}
	\newcommand{\ees}{\end{eqnarray*}}
\newcommand{\bpm}{\begin{pmatrix}}
	\newcommand{\epm}{\end{pmatrix}}
	\newcommand{\bbm}{\begin{bmatrix}}
	\newcommand{\ebm}{\end{bmatrix}}

\def\diag{{\rm diag}\,}

\begin{document}

%% End-Of-Header

\title[Switching and partially switching the hypercube]{Switching and partially switching the hypercube while maintaining perfect state transfer}

\author[S.\ Kirkland]{Steve Kirkland$^1$}
\author[S.\ Plosker]{Sarah Plosker$^{2,1}$}
\author[X.~Zhang]{Xiaohong Zhang$^{1}$}

\address{$^1$Department of Mathematics, University of Manitoba, Winnipeg, MB, Canada  R3T 2N2}
\address{$^2$Department of Mathematics \& Computer Science, Brandon University, Brandon,
 MB, Canada R7A 6A9}


\begin{abstract} A graph is said to exhibit perfect state transfer (PST) if one of its corresponding Hamiltonian matrices, which are based on the vertex-edge structure of  the graph, gives rise to PST in a quantum information-theoretic context, namely with respect to inter-qubit interactions of a quantum  system. We perform various perturbations to the hypercube graph---a graph that is known to exhibit PST---to create graphs that maintain many of the same properties of the hypercube, including PST as well as the distance for which PST occurs. We show that the sensitivity with respect to readout time errors remains unaffected for the vertices involved in PST. We give motivation for when these perturbations may be physically desirable or even necessary.
\end{abstract}
\maketitle
\section{Introduction}        
 Undirected connected graphs are used as models for quantum spin networks, and in particular to model  inter-qubit interactions of quantum registers and processors within a quantum computer. The transfer of quantum states from one location to another within a quantum computer is then analyzed by way of a Hamiltonian $\mathcal{H}$, which is a   matrix describing the total
amount of energy of a quantum system, and, depending on the dynamics of the quantum system, is typically taken to be either the adjacency matrix or Laplacian matrix associated to the graph. 

%Here, we focus on dynamics that result in the consideration of the associated Laplacian matrix. 
Recent work explores hypercubes (also called $n$-cubes) \cite{Cubes} and the more general notion of cubelike graphs \cite{Cubelike,  MBHadamards} as a means of achieving perfect state transfer (PST): a quantum state placed at a particular vertex of the spin network is transmitted perfectly (up to a global phase) to another vertex at time $t=t_0$.  In \cite{MBHadamards}, the authors  make  use of results in \cite{BFK} involving    Laplacians that can be diagonalized by a Hadamard matrix to create a variety of new graphs having PST; again, particular attention is paid to the hypercube and, more generally, to cubelike graphs. One is often interested in sending the state as far as possible along the spin network (maximizing the distance between the  vertices of the sender and receiver), and ensuring that the state transfer is as insensitive as possible to errors in the readout time (that is, if one has PST at time $t=t_0$, then at time $t=t_0\pm \epsilon$, one would hope to have near-perfect state transfer, for small epsilon).

We present new results on variants of the hypercube while avoiding much of the heavy machinery from graph theory. We use a graph operation called Godsil-McKay (GM) switching   that perturbs a graph by removing and creating edges based on certain criteria that a partition of the vertex set must satisfy. We apply GM switching to the hypercube as a means of constructing non-isomorphic graphs that have many of the same nice properties of the hypercube, including PST. 

Other recent work   perturbs the graph, often taken to be a quantum spin chain modeled by a path,  by adding loops (corresponding to energy potentials) at certain vertices \cite{potentials, KemptonPST, RV}. This type of perturbation has garnered much attention as a powerful tool to improve quantum state transfer capabilities of quantum networks.  Other work (\cite{CDDEKL, signed} and others) focuses more on manufacturing coupling strengths (corresponding to the weights of the edges of the graph) to achieve perfect state transfer. Small changes in the edge weights lead to manufacturing errors; an analysis of how such errors decrease the probability of quantum state transfer can be found in \cite{Steve2015}. Other more applied literature discusses enhancing the probability of state transfer by way of partially collapsing measurements, weak measurement strength, and quantum measurement reversal \cite{BA}, radio frequency pulses in NMR \cite{CRC}, quantum error correction \cite{Kay15}, among other useful perturbations, both at the local (individual vertices) and global (quantum system) levels. 

In Section \ref{sec:prelim}, we give the necessary graph theory and linear algebra background for this work. In Section \ref{sec:GMswitching}, we use Godsil-McKay switching to construct a graph (the {\it{switched $n$-cube}}) of order $2^n$ for $n>4$ that has many of the same properties as the $n$-cube (notably, it is cospectral to the $n$-cube and exhibits PST, with distance $n$ between PST pairs), but is nevertheless not isomorphic to the $n$-cube, and is not Hadamard-diagonalizable, unlike the $n$-cube. In Section \ref{sec:partialS}, we then consider partially switched $n$-cubes, which generalize the process of GM switching on the $n$-cube by considering 
it as the Cartesian product of the $(n-4)$-cube with the 4-cube, and performing GM switching on some copies of the 4-cube. %the corresponding adjacency matrix in block form and performing GM switching to only some of the blocks. 
These new graphs are not cospectral to the $n$-cube in general, but do exhibit PST (though in significantly fewer pairs of vertices).  We further generalize this by replacing each copy of the 4-cube by a convex combination of the 4-cube and the switched 4-cube; %diagonal block in the adjacency matrix with convex combinations of the $4$-cube and a partially switched $4$-cube adjacency matrices; 
we also generalize it to a time-dependent Hamiltonian (see, e.g.\ \cite{Pan}) that alternates between the various graphs considered.  We give motivation as to why these families of graphs might be useful in practice; in particular, we conduct a sensitivity analysis with respect to readout time errors in Section \ref{sec:analysis}.

\section{Preliminaries}\label{sec:prelim}
We consider only   unweighted (with the exception of a generalization in Section \ref{sec:partialS}), undirected, simple, connected graphs herein. %Our work is based on the  $n$-cube, which can be constructed  by labelling its $2^n$ vertices $v=1,2,\dots, 2^n$ by the binary representations of the numbers $v-1$ respectively, and then connecting vertices with an edge precisely when their Hamming distance is 1. The $n$-cube is $n$-regular:  each vertex has exactly $n$ adjacent vertices. The graph of the $n$-cube is denoted $Q_n$. 
Given a graph $G$ on $m$ vertices, its corresponding \textit{adjacency matrix} is an $m \times m$ matrix $A(G)=[a_{jk}]$ with $a_{jk}=1$ if vertices $j$ and $k$ are  adjacent, and $a_{jk}=0$ otherwise (in general, $a_{jk}$ represents the weight of the edge between vertices $j$ and $k$ in a weighted graph).   The Laplacian matrix $L(G)$ corresponding to  $G$ is defined as the $m\times m$ matrix
$L(G)=D(G)-A(G)$, where $D(G)$ (the degree matrix) is the diagonal matrix of row sums of $A(G)$.  Since the hypercube is regular,  many of the properties of the Laplacian matrix $L(G)$   are shared by the adjacency matrix $A(G)$, and we focus on the latter  herein. In particular, PST with respect to the Laplacian matrix will occur if and only if PST with respect to the adjacency matrix occurs, for the same vertex pairs; this statement is in general not true if regularity is dropped.

The graph $G$ is used to represent a single-excitation spin network, with each electron being represented by a vertex of the graph, and coupling strengths being represented by the weight of the edge between the two interacting electrons (adjacent vertices). The system is given by the Hilbert space ${\mathbb{C}^{2}}^{\otimes n}$; here, we are interested in the $XY$-interaction model (alternative terminology is that the spin network has $XX$ couplings). The total energy of the system is thus given by the Hamiltonian
\[
\mathcal{H}=\frac12\sum_{\{j,k\}\in E(G)} w_{jk}(X_jX_k+Y_jY_k)
\]
where $\{j,k\}\in E(G)$ means that there is an edge between vertices $j$ and $k$ in the graph, having edge weight $w_{jk}$ (in nearly all of what follows, $w_{jk}=1$), and $X_j$ and $Y_j$ are the standard Pauli matrices acting on the $j$-th copy of $\mathbb{C}^2$. Although in the above, $\mathcal H$ is a matrix of order $2^n$, since we are considering only the single excitation subspace $\mathbb{C}^n$ spanned by the standard basis $\{\ket{1}, \dots, \ket{n}\}\subset \mathbb{C}^n$, the state transfer dynamics are completely determined by the evolution within this $n$-dimensional subspace, and the Hamiltonian can be represented by  $\mathcal{H}=A(G)$ when considering XY interactions, as per  \cite{Hamiltonian}.  

The \emph{probability (or fidelity) of state transfer} is %or \emph{fidelity} between vertices $j$ and $k$ is defined as 
\begin{equation}\label{eq:fidelity}
p_{j,k}(t)= |\bra{j}e^{it\mathcal{H}}\ket{k}|^2=|\bra{k}e^{it\mathcal{H}}\ket{j}|^2,
\end{equation}
 where $\mathcal{H}$ is the Hamiltonian of the system (since we are considering XY interactions, $\mathcal{H}=A(G)$).  Since $\mathcal{H}$ is symmetric, we can consider either the $(j,k)$ or $(k,j)$ entry of  the unitary matrix $U(t)=e^{it\mathcal{H}}$. If there exists a time $t=t_0$ for which $p_{j,k}(t_0)=1$, then we say that the vertices $j$ and $k$ exhibit PST (or that the graph has PST, or that $\{j,k\}$ is a PST vertex pair, or that $j$ and $k$ pair up  
and have PST). If the graph theoretic distance between vertices $j$ and $k$ is $\ell$ (that is, the  minimum  number of edges in a   path joining $j$ and $k$), then we say that the PST distance is $\ell$.  %In this paper we focus on regular graphs, so the Laplacian matrix $L$ equals to $dI-A$, where $d$ is the degree of the regular graph, $A$ is its adjacency matrix, and $I$ is the identity matrix. 
For regular graphs, which are our focus here,  a graph exhibits PST with respect to the Laplacian matrix if and only if it exhibits PST with respect to the adjacency matrix, so focussing on the adjacency matrix is not a restriction in this setting.

We use $I_n$ to denote the identity matrix of size $2^n$,   and $\mathbf{1}_{m}$ to denote the unnormalized maximally mixed state (all-ones vector) of length $2^m$. 
A \emph{Hadamard matrix} of order $m$ is an $m\times m$ square  $(1,-1)$ matrix whose columns are pairwise orthogonal. The \emph{standard Hadamard matrices} of order $2^n$ for $n\in \mathbb{N}$ are defined recursively:
 let
 \[
 H_1=\begin{bmatrix}
 1 & 1\\
 1 & -1
 \end{bmatrix},
\quad
 H_2=\begin{bmatrix}
 H_1&H_1\\H_1&-H_1
 \end{bmatrix}=\begin{bmatrix}
  1 &  1 & 1 &  1\\
  1 & -1 & 1 & -1\\
  1 &  1 &-1 & -1\\
  1 & -1 &-1 & 1\\
  \end{bmatrix},
  \]
  and then
  \[
  H_n=\begin{bmatrix}
  H_{n-1} & H_{n-1}\\
  H_{n-1} &-H_{n-1}\\
  \end{bmatrix}=H_1\otimes H_{n-1}=H_1^{\otimes n},
  \]
  for $2 \leq n\in\mathbb{N}$, where $H_1^{\otimes n}$ denotes the Kronecker product of $H_1$ with itself $n$ times.

A graph $G$ on $m$ vertices is \emph{Hadamard-diagonalizable} if we can write $L(G)=\frac1m H\Lambda H^T$, where $\Lambda$ is a diagonal matrix of eigenvalues of $L(G)$, and $H$ is a Hadamard matrix.  
From \cite{BFK} we know that any Hadamard diagonalizable graph is regular, so a graph is Hadamard diagonalizable if and only if its adjacency matrix is diagonalizable by some Hadamard matrix. In this paper, we will make use of the adjacency matrix of a graph to check whether it is Hadamard diagonalizable or not.
Here, we focus on standard Hadamard matrices, and so $m=2^n$ for some $n\in \mathbb{N}$.

A useful family of Hadamard diagonalizable graphs is the
family of cubelike graphs \cite{BGS08}: Take a set
$C\subset \mathbb{Z}_{2}^{n}=\mathbb{Z}_{2}\times \cdots \times
\mathbb{Z}_{2}$ ($n$ times), where $C$ does not contain the all-zeros
vector. Construct the \emph{cubelike graph} $G(C)$ with vertex set
$V=\mathbb{Z}_{2}^{n}$ and two elements of $V$ are adjacent if and only
if their difference is in $C$. The set $C$ is called the
\emph{connection set} of the graph $G(C)$. From the definition we can see that $G(C)$ is a  $|C|$-regular graph. A cubelike graph is connected if and only if its connection set $C$ contains a basis of                            $\mathbb{Z}_{2}^{n}$ when viewed as a vector space or $C$ generates $\mathbb{Z}_2^n$ when viewed as a group  \cite[Ch.\ 3]{C.Galgg}.
Since  $\mathbb{Z}_{2}^{n}$ is $n$-dimensional, we know a connected cubelike graph on $2^n$ vertices is regular with degree at least $n$.

An unweighted graph $G$ is diagonalizable by the standard Hadamard matrix if and only if $G$ is a cubelike graph  \cite{MBHadamards}. %the Laplacian matrix $L$ of                     a. 
Combined with the above information, we know that for any positive integer $n$, no connected graphs on $2^n$ vertices that have fewer edges than the $n$-cube are diagonalizable by the standard Hadamard matrix $H_n$.
%no connected graphs which are strictly sparser than hypercubes are diagonalizable by the standard Hadamard matrices.
So we cannot perturb the hypercube by deleting edges, without adding edges as well, while still maintaining PST, connectivity, and diagonalizability by the standard Hadamard. We summarize this in the following proposition.

\vspace*{12pt}
\noindent
{\bf Proposition~1:}  
No connected proper subgraphs of hypercubes are diagonalizable by the standard Hadamard matrix $H_n$.

\vspace*{12pt}

The switched cube discussed in Section \ref{sec:GMswitching} maintains the same sparsity structure (the same number of edges)  as the $n$-cube $Q_n$.

\section{A PST Graph Cospectral to the $n$-cube}\label{sec:GMswitching}

 In this section, we give a cospectral mate of the $n$-cube that is no longer Hadamard diagonalizable, but is $n$-regular, exhibits PST, and has PST distance $n$.   

The Cartesian product of two graphs $G_1$ and $G_2$ gives a new graph $G_1\square G_2$ with vertex set $V(G_1)\times V(G_2)$ and two vertices $(j_1,j_2)$ and $(k_1, k_2)$ are adjacent in $G_1\square G_2$ provided either 
   $j_1=k_1$ and $j_2$ is adjacent to $k_2$ in $G_2$, or 
    $j_2=k_2$ and $j_1$ is adjacent to $k_1$ in $G_1$.
    
We will make use of Godsil-McKay (GM) switching  
  \cite{GM}:  
Let $G$ be a graph and let $\pi=(C_1,C_2,\cdots,C_k,D)$ be a partition of the vertex set $V(G)$.
Suppose that, whenever $1\leq i,j\leq k$ and $v\in D$, we have: \\
(a) any two vertices in $C_i$ have the same number of neighbours in $C_j$, and\\
(b) $v$ has either 0, $n_i/2$ or $n_i$ neighbours in $C_i$, where $n_i=|C_i|$.\\
The graph $G^{(\pi)}$ formed by local switching in $G$ with respect to $\pi$ is obtained from $G$ as follows: for each $v\in D$ and $1\leq i\leq k$ such that $v$ has $n_i/2$ neighbours in $C_i$, delete these $n_i/2$
edges and join $v$ instead to the other $n_i/2$ vertices in $C_i$. 
The graphs $G^{(\pi)}$ and $G$ are cospectral.

\begin{figure} [h] %[htbp]
\vspace*{13pt}
 \begin{center}
 \definecolor{xdxdff}{rgb}{0.49019607843137253,0.49019607843137253,1.}
\definecolor{qqqqff}{rgb}{0.,0.,1.}
\begin{tikzpicture}[line cap=round,line join=round,>=triangle 45,x=1.0cm,y=1.0cm,scale=0.45]
\clip(1.5,-10) rectangle (48,9);
\draw (7.000815999999996,5.943981999999998)-- (4.,4.);
\draw (7.000815999999996,5.943981999999998)-- (6.,4.);
\draw (7.000815999999996,5.943981999999998)-- (8.,4.);
\draw (7.000815999999996,5.943981999999998)-- (10.,4.);
\draw (4.,4.)-- (2.,0.);
\draw (4.,4.)-- (4.,0.);
\draw (4.,4.)-- (10.,0.);
\draw (6.,4.)-- (2.,0.);
\draw (6.,4.)-- (6.,0.);
\draw (6.,4.)-- (12.,0.);
\draw (8.,4.)-- (4.,0.);
\draw (8.,4.)-- (8.,0.);
\draw (8.,4.)-- (12.,0.);
\draw (10.,4.)-- (6.,0.);
\draw (10.,4.)-- (8.,0.);
\draw (10.,4.)-- (10.,0.);
\draw (2.,0.)-- (8.,-4.);
\draw (2.,0.)-- (10.,-4.);
\draw (4.,0.)-- (6.,-4.);
\draw (4.,0.)-- (10.,-4.);
\draw (6.,0.)-- (4.,-4.);
\draw (6.,0.)-- (8.,-4.);
\draw (8.,0.)-- (4.,-4.);
\draw (8.,0.)-- (6.,-4.);
\draw (10.,0.)-- (6.,-4.);
\draw (10.,0.)-- (8.,-4.);
\draw (12.,0.)-- (4.,-4.);
\draw (12.,0.)-- (10.,-4.);
\draw (4.,-4.)-- (7.000815999999996,-5.8273819999999965);
\draw (6.,-4.)-- (7.000815999999996,-5.8273819999999965);
\draw (8.,-4.)-- (7.000815999999996,-5.8273819999999965);
\draw (10.,-4.)-- (7.000815999999996,-5.8273819999999965);
\draw (21.173304,5.9147)-- (18.172488,3.970717999999999);
\draw (21.173304,5.9147)-- (20.172488,3.970717999999999);
\draw (21.173304,5.9147)-- (22.172488,3.970717999999999);
\draw (21.173304,5.9147)-- (24.,4.);
\draw (16.172488,-0.029282000000000585)-- (22.172488,-4.029282);
\draw (16.172488,-0.029282000000000585)-- (24.172488,-4.029282);
\draw (18.,0.)-- (20.172488,-4.029282);
\draw (18.,0.)-- (24.172488,-4.029282);
\draw (20.172488,-0.029282000000000585)-- (18.172488,-4.029282);
\draw (20.172488,-0.029282000000000585)-- (22.172488,-4.029282);
\draw (22.172488,-0.029282000000000585)-- (18.172488,-4.029282);
\draw (22.172488,-0.029282000000000585)-- (20.172488,-4.029282);
\draw (24.,0.)-- (20.172488,-4.029282);
\draw (24.,0.)-- (22.172488,-4.029282);
\draw (18.172488,-4.029282)-- (21.173304,-5.856663999999997);
\draw (20.172488,-4.029282)-- (21.173304,-5.856663999999997);
\draw (22.172488,-4.029282)-- (21.173304,-5.856663999999997);
\draw (24.172488,-4.029282)-- (21.173304,-5.856663999999997);
\draw (26.,0.)-- (18.172488,-4.029282);
\draw (26.,0.)-- (24.172488,-4.029282);
\draw (18.172488,3.970717999999999)-- (20.172488,-0.029282000000000585);
\draw (18.172488,3.970717999999999)-- (22.172488,-0.029282000000000585);
\draw (18.172488,3.970717999999999)-- (26.,0.);
\draw (20.172488,3.970717999999999)-- (18.,0.);
\draw (20.172488,3.970717999999999)-- (22.172488,-0.029282000000000585);
\draw (20.172488,3.970717999999999)-- (24.,0.);
\draw (22.172488,3.970717999999999)-- (16.172488,-0.029282000000000585);
\draw (22.172488,3.970717999999999)-- (20.172488,-0.029282000000000585);
\draw (22.172488,3.970717999999999)-- (24.,0.);
\draw (24.,4.)-- (26.,0.);
\draw (24.,4.)-- (18.,0.);
\draw (24.,4.)-- (16.172488,-0.029282000000000585);
%\draw (5.5,-7) node[anchor=north west] {\small{4-cube $Q_4$}};
%\draw (17.2,-7) node[anchor=north west] {\small{Switched 4-cube $\tilde{Q}_4=Q_4^{(\pi)}$}};
\begin{scriptsize}
\draw [fill=qqqqff] (7.000815999999996,5.943981999999998) circle (2.5pt);
\draw[color=qqqqff] (7.4693279999999955,6.471058000000004) node {$1$};
\draw [fill=qqqqff] (4.,4.) circle (2.5pt);
\draw[color=qqqqff] (4.541127999999996,4.626292000000005) node {$2$};
\draw [fill=qqqqff] (6.,4.) circle (2.5pt);
\draw[color=qqqqff] (6.47374,4.538446000000005) node {$3$};
\draw [fill=qqqqff] (8.,4.) circle (2.5pt);
\draw[color=qqqqff] (8.464915999999995,4.538446000000005) node {$4$};
\draw [fill=qqqqff] (10.,4.) circle (2.5pt);
\draw[color=qqqqff] (10.456091999999995,4.538446000000005) node {$5$};
\draw [fill=qqqqff] (2.,0.) circle (2.5pt);
\draw[color=qqqqff] (2.4621059999999964,0.5268120000000058) node {$6$};
\draw [fill=qqqqff] (4.,0.) circle (2.5pt);
\draw[color=qqqqff] (4.482563999999996,0.5268120000000058) node {$7$};
\draw [fill=qqqqff] (6.,0.) circle (2.5pt);
\draw[color=qqqqff] (6.47374,0.5268120000000058) node {$10$};
\draw [fill=qqqqff] (8.,0.) circle (2.5pt);
\draw[color=qqqqff] (8.464915999999995,0.5268120000000058) node {$11$};
\draw [fill=qqqqff] (10.,0.) circle (2.5pt);
\draw[color=qqqqff] (10.456091999999995,0.5268120000000058) node {$9$};
\draw [fill=qqqqff] (12.,0.) circle (2.5pt);
\draw[color=qqqqff] (12.535113999999995,0.6146580000000057) node {$8$};
\draw [fill=qqqqff] (4.,-4.) circle (2.5pt);
\draw[color=qqqqff] (4.482563999999996,-3.4848219999999928) node {$15$};
\draw [fill=qqqqff] (6.,-4.) circle (2.5pt);
\draw[color=qqqqff] (6.47374,-3.4848219999999928) node {$14$};
\draw [fill=qqqqff] (8.,-4.) circle (2.5pt);
\draw[color=qqqqff] (8.464915999999995,-3.4848219999999928) node {$13$};
\draw [fill=qqqqff] (10.,-4.) circle (2.5pt);
\draw[color=qqqqff] (10.456091999999995,-3.4848219999999928) node {$12$};
\draw [fill=qqqqff] (7.000815999999996,-5.8273819999999965) circle (2.5pt);
\draw[color=qqqqff] (7.4693279999999955,-5.300305999999992) node {$16$};
\draw [fill=qqqqff] (21.173304,5.9147) circle (2.5pt);
\draw[color=qqqqff] (21.70038,6.529622000000003) node {$1$};
\draw [fill=qqqqff] (18.172488,3.970717999999999) circle (2.5pt);
\draw[color=qqqqff] (18.71361599999999,4.59701) node {$2$};
\draw [fill=qqqqff] (20.172488,3.970717999999999) circle (2.5pt);
\draw[color=qqqqff] (20.70479199999999,4.59701) node {$3$};
\draw [fill=qqqqff] (22.172488,3.970717999999999) circle (2.5pt);
\draw[color=qqqqff] (22.69596799999999,4.59701) node {$4$};
\draw [fill=qqqqff] (24.,4.) circle (2.5pt);
\draw[color=qqqqff] (24.54073399999999,4.626292000000005) node {$5$};
\draw [fill=qqqqff] (16.172488,-0.029282000000000585) circle (2.5pt);
\draw[color=qqqqff] (16.69315799999999,0.5853760000000058) node {$6$};
\draw [fill=qqqqff] (18.,0.) circle (2.5pt);
\draw[color=qqqqff] (18.53792399999999,0.6146580000000057) node {$7$};
\draw [fill=qqqqff] (20.172488,-0.029282000000000585) circle (2.5pt);
\draw[color=qqqqff] (20.70479199999999,0.5853760000000058) node {$10$};
\draw [fill=qqqqff] (22.172488,-0.029282000000000585) circle (2.5pt);
\draw[color=qqqqff] (22.69596799999999,0.5853760000000058) node {$11$};
\draw [fill=qqqqff] (24.,0.) circle (2.5pt);
\draw[color=qqqqff] (24.54073399999999,0.6146580000000057) node {$9$};
\draw [fill=qqqqff] (18.172488,-4.029282) circle (2.5pt);
\draw[color=qqqqff] (18.71361599999999,-3.4262579999999927) node {$15$};
\draw [fill=qqqqff] (20.172488,-4.029282) circle (2.5pt);
\draw[color=qqqqff] (20.70479199999999,-3.4262579999999927) node {$14$};
\draw [fill=qqqqff] (22.172488,-4.029282) circle (2.5pt);
\draw[color=qqqqff] (22.69596799999999,-3.4262579999999927) node {$13$};
\draw [fill=qqqqff] (24.172488,-4.029282) circle (2.5pt);
\draw[color=qqqqff] (24.68714399999999,-3.4262579999999927) node {$12$};
\draw [fill=qqqqff] (21.173304,-5.856663999999997) circle (2.5pt);
\draw[color=qqqqff] (21.70038,-5.2417419999999915) node {$16$};
\draw [fill=qqqqff] (26.,0.) circle (2.5pt);
\draw[color=qqqqff] (26.47334599999999,0.5268120000000058) node {$8$};
\end{scriptsize}
\end{tikzpicture}
 \end{center}
\vspace*{13pt}
\caption{ \label{4-cubepar}The 4-cube $Q_4$ (left) and the switched 4-cube $\tilde{Q}_4$ (right).}
\end{figure}
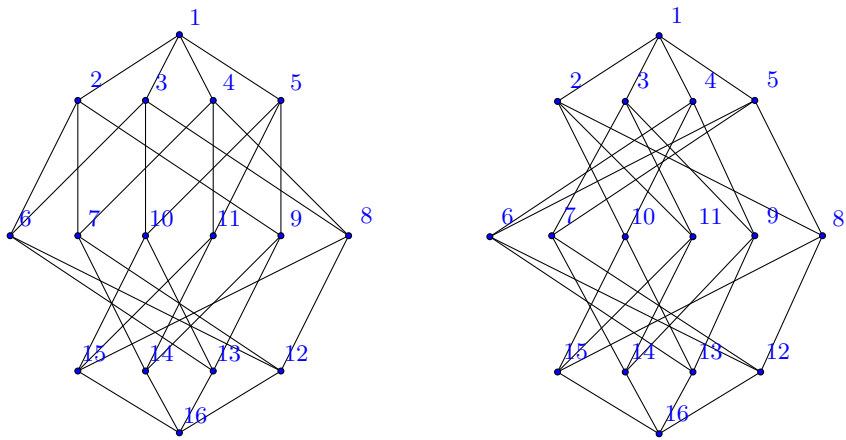

In \cite[Section 1.8]{ABWH}, a cospectral mate of the 4-cube is given, and can be seen to be an example   of Godsil-McKay (GM) switching. With labeling as in Fig.~\ref{4-cubepar}, the partition $\pi$  has  $C_1=\{1\}$, $C_2=\{6,7,8,9,10,11\}$, $C_3=\{12,13,14,15\}$, $C_4=\{16\}$, and $D=\{2,3,4,5\}$; denote the graph  $Q_4^{(\pi)}$ by $\tilde{Q}_4$, and call it the switched 4-cube.
Note that $Q_n=Q_{n-4}\square Q_4$, that is, the $n$-cube can be seen as $2^{n-4}$ copies of the 4-cube connected in a specific way (in fact, according to the $(n-4)$-cube). If we partition each of the $2^{n-4}$ copies of 4-cube according to the above partition $\pi$, we get a partition $\pi_1$ of $V(Q_n)$, and it is an equitable partition. Now taking the union of all the $2^{n-4}$ copies of the $D$ cell and keeping all the other cells unchanged, we get a new partition $\pi_2$ of $V(Q_n)$, which satisfies the GM switching conditions. Denote $Q_n^{(\pi_2)}$ as $\tilde{Q}_n$  and call it the switched $n$-cube. From the construction we know that $\tilde{Q}_n=Q_{n-4}\square \tilde{Q_4}$. For this switched $n$-cube, we order the vertices of $\tilde{Q}_4$ as  in Fig.~\ref{4-cubepar}, order the vertices of $Q_{n-4}$ in increasing order of their binary representations, and finally order the vertices of the Cartesian product $\tilde{Q}_n=Q_{n-4}\square \tilde{Q_4}$ in accordance with the dictionary ordering, that is, $A(\tilde{Q}_n)=A(Q_{n-4})\otimes I_4+I_{n-4}\otimes A(\tilde{Q}_4)$ \cite{Fiedleral}. Order the vertices of the $n$-cube accordingly.

\iffalse
since we only make use of this partition, we use the labeling $\tilde{Q}_4$ for the graph $Q_4^{(\pi)}$ formed by local switching in the 4-cube $Q_4$. We will build switched $n$-cubes $\tilde{Q}_n$ from the switched 4-cube $\tilde{Q}_4$ via the Cartesian product construction $Q_n=Q_{n-4}\square Q_4$, from which we obtain  $\tilde{Q}_n=Q_{n-4}\square \tilde{Q}_4$, although we note that, equivalently, $\tilde{Q}_n$ can  be obtained from $Q_n$ through GM switching using a particular partition $\pi$. %The switched $n$-cube $\tilde{Q}_n$ is a cospectral mate of the $n$-cube $Q_n$, and like the switched $4$-cube, it can be obtained from the $n$-cube through GM switching, which relies on a partition $\pi=(V_1, V_2, \dots, V_k, D)$ of $V(G)$; we simply note that $D$ for the $n$-cube is the union of all $2^{n-4}$ copies of $D$ for the 4-cube, and all $V_j$ remain unchanged.   
 \fi

 Denote the   adjacency matrix of $Q_n$ by $C_n$, and the adjacency matrix of $\tilde{Q}_n$ by $\tilde{C}_n$. 
The non-isomorphism of $\tilde{Q}_n$ and $Q_n$ can be seen directly from the fact that they exhibit different PST properties: namely, 
they have different numbers of PST vertex pairs.

\vspace*{12pt}
\noindent
{\bf Theorem~2:} \label{thm:switchedPSTpairs}
For $n\geq 4$, exactly half of the vertices of the switched $n$-cube $\tilde{Q}_n$ pair up and have PST between each other at time $\pi/2$. 

\vspace*{12pt}
\noindent
{\bf Proof:} 
Since the eigendecomposition of $\tilde{C}_4$ is known in closed form, we may explicitly compute $e^{i\tilde{C}_4 \pi/2}$, also in closed form. From that explicit computation, we can see that there is PST between vertices 1 and 16,  6 and 11, 7 and 10, 8 and 9 in $\tilde{Q}_4$ at time $\pi/2$ (whereas the $4$-cube has PST between vertices $j$ and $17-j$ for each $j$ in this ordering); exactly half (8 out of 16) of the vertices pair up. 
%Since $Q_{n-4}\square \tilde{Q}_4$, we know the adjacency matrix $\tilde{C}_{n}$ of $\tilde{Q}_{n}$ is
Recall $\tilde{C}_{n}=C_{n-4}\otimes I_4+I_{n-4}\otimes\tilde{C}_4 $.
Therefore 
\begin{eqnarray*}
U(\pi/2)&=&e^{i (\pi/2)\tilde{C}_{n}}=e^{i(\pi/2)(C_{n-4}\otimes I_4+I_{n-4}\otimes\tilde{C}_4 )}\\
&=&e^{i(\pi/2)C_{n-4}\otimes I_4}e^{i(\pi/2)I_{n-4}\otimes\tilde{C}_4}=e^{i(\pi/2)C_{n-4}}\otimes e^{i(\pi/2)\tilde{C}_4}\\
&=&(i)^{(n-4)}\begin{bmatrix}0& 0&\cdots &0 &e^{i(\pi/2)\tilde{C}_4}\\
0& 0&\cdots  &e^{i(\pi/2)\tilde{C}_4}&0\\
\cdots &\cdots &\cdots &\cdots &\cdots\\
 0&e^{i(\pi/2)\tilde{C}_4}&\cdots  &0 &0\\
e^{i(\pi/2)\tilde{C}_4}&0&\cdots  &0&0\\
\end{bmatrix},
\end{eqnarray*}
where the third equality is based on the fact that $\tilde{C}_k\otimes I_1$ and $I_k\otimes C_1$ commute.
Hence if half  the rows of $e^{i(\pi/2)\tilde{C}_4}$ have an offdiagonal entry with modulus 1, then so does $e^{i(\pi/2)\tilde{C}_n}$, i.e., half of the vertices of $\tilde{Q}_n$ pair up and have PST between each other at time $\pi/2$.
$\square$

From \cite{BGS08,Cubelike} we know that all vertices of the $n$-cube pair up to have PST at time $\pi/2$ (that is, each vertex of the $n$-cube is part of a vertex pair for which PST occurs at time $\pi/2$: namely, PST occurs between vertex $k$ and vertex $2^n+1-k$ for $k=1,\ldots, 2^n$, where again the vertices of the $n$-cube and switched $n$-cube are ordered as mentioned above Theorem~\ref{thm:switchedPSTpairs}. We thus have the following corollary, implying that the two graphs $Q_n$ and $\tilde{Q}_n$ are  non-isomorphic.

\vspace*{12pt}
\noindent
{\bf Corollary~3:}
There are half as many  vertex pairs for which PST occurs (at time $\pi/2$) for the switched $n$-cube as there are for the $n$-cube.
\vspace*{12pt}
\noindent

\section{Partial Switching and PST}\label{sec:partialS}

%In Section \ref{sec:GMswitching}, we explored the application of GM switching to the $n$-cube to produce a switched $n$-cube.  
In Section \ref{sec:GMswitching}, we explored the PST property of the cospectral mate $\tilde{Q}_n$ of the $n$-cube, which can be obtained from each other through GM switching. In this Section, we continue considering the Cartesian product construction of the $n$-cube ($Q_n=Q_{n-4}\square Q_4$) so that the corresponding adjacency matrix is seen to be a block matrix, with each block of size $16\times16$, and the diagonal blocks are copies of the adjacency matrix $C_4$ of $Q_4$. We then perform GM switching to some (but not all) copies of $Q_4$ inside the $n$-cube, i.e., some diagonal blocks of $C_n$ are changed from $C_4$ to $\tilde{C}_4$. We call this \emph{partial switching}, and we analyze the property of PST for these partially switched $n$-cubes.

\subsection{Construction} Let $n\geq 4$. 
Let $A_{n,1}=C_n$ be the adjacency matrix of the $n$-cube, and let $A_{n,2}=\tilde{C}_n$ be the adjacency matrix of the switched $n$-cube. %, with
The adjacency matrix of the $n$-cube is $A_{n,1}=A_{n-4,1}\otimes I_4+I_{n-4}\otimes A_{4,1}=\diag(A_{4,1},\cdots,$ $ A_{4,1})+A_{n-4,1}\otimes I_4$.
For example, $A_{5,1}=\begin{bmatrix}A_{4,1} &I_4\\I_4 &A_{4,1}\end{bmatrix}$, and 
$A_{6,1}=$ $\begin{bmatrix}A_{4,1}&I_4&I_4&0\\ I_4&A_{4,1}&0&I_4\\I_4&0&A_{4,1}&I_4\\0&I_4&I_4&A_{4,1}\end{bmatrix}$. 
For the switched $n$-cube, we just need to replace every occurrence of $A_{4,1}$ with $A_{4,2}$.
%Here we also consider partially switched $n$-cubes: 
Now, if we replace some of the diagonal blocks $A_{4,1}$ in $A_{n,1}$
by $A_{4,2}$, we get the adjacency matrix of a partially switched $n$-cube. For $n=5$, by a simple reordering of the copies of  $A_{4,1}$ and $A_{4,2}$, it is clear that  $\begin{bmatrix}A_{4,1} &I_4\\I_4 &A_{4,2}\end{bmatrix}$ and $\begin{bmatrix}A_{4,2} &I_4\\I_4 &A_{4,1}\end{bmatrix}$ are isomorphic; but they are not isomorphic to the 5-cube or the switched 5-cube (by checking that they have different spectrum or by the result of Example~7 in Section \ref{sec:wvmPST} below we know they have fewer vertices that exhibit PST).

For $n=6$, there is a unique (up to isomorphism, which can be accomplished by  reordering the copies of $A_{4,1}$ and $A_{4,2}$) partially switched 6-cube with exactly one copy of the 4-cube, say $A_{6,3}=\diag(A_{4,1},A_{4,2},A_{4,2},A_{4,2})+A_{2,1}\otimes I_4$ and a unique (again, up to isomorphism) partially switched 6-cube with three copies of the 4-cube, say  $A_{6,4}=\diag(A_{4,2},A_{4,1},A_{4,1},A_{4,1})+A_{2,1}\otimes I_4$; furthermore, $G(A_{6,4})$ can be obtained from $G(A_{6,3})$ by performing the GM switching. When there are two copies of $A_{4,1}$ and two copies of $A_{4,2}$, there are two nonisomorphic partially switched 6-cubes: $A_{6,5}=\diag(A_{4,1},A_{4,1},A_{4,2},A_{4,2})+A_{2,1}\otimes I_4$, 
$A_{6,6}=\diag(A_{4,1},A_{4,2},A_{4,2},$ $A_{4,1})+A_{2,1}\otimes I_4$, and these two graphs are not even cospectral. For $n\geq 7$, there are more different types of partially switched $n$-cubes.

%This holds true for $n=6$ as well: any two partially switched 6-cubes with the same number of occurrence of $A_{4,1}$ (and therefore the same number of occurrences of $A_{4,2}$)  are isomorphic to each other simply by permuting the blocks (this can be easily checked as the dimension is relatively small).

%We consider  arbitrary partially switched $n$-cubes with adjacency matrix of the form 
%$\bbm A_{(n-2)*}&I&I&0\\ I&A_{(n-2)*}&0&I\\I&0&A_{(n-2)*}&I\\0&I&I&A_{(n-2)*}\ebm$ where $A_{(n-2)*}$ is the adjacency matrix of the  $(n-2)$-cube or of the switched $(n-2)$-cube or of some partially switched $(n-2)$-cube. Note that if all blocks are $A_{(n-2)1}$ then this is just the adjacency matrix of the $n$-cube; if all blocks are  $A_{(n-2)2}$ then this is just the adjacency matrix of the switched $n$-cube. It turns out that, for $n\geq 7$, partially switched $n$-cubes
%are not uniquely determined (up to isomorphism) by the number of appearance of $A_{4,1}$ anymore. For example, the partially switched 7-cube with adjacency matrix $\diag(A_{4,1},A_{4,2},A_{4,1},A_{4,1},A_{4,1},A_{4,1},A_{4,2},A_{4,1})+A_{31}\otimes I_4$ and  the partially switched 7-cube with adjacency matrix       $\diag(A_{4,1},A_{4,1},A_{4,1},A_{4,1},A_{4,1},A_{4,1},A_{4,2},A_{4,2})+A_{31}\otimes I_4$ are not isomorphic to each other.

The partially switched $n$-cubes are no longer cospectral to the $n$-cube in general, but still exhibit PST (though the number of PST vertex pairs is significantly fewer  than in the $n$-cube), and are not cubelike graphs, for the following reason.  It can be shown that there is no isomorphism of a partially switched $n$-cube that maps vertex 16   to vertex 6 within any copy of $\tilde{Q}_4$, which shows that the graph is not vertex-transitive, and hence it is not a cubelike graph. % (since any Cayley graph is vertex-transitive).

\subsection{Which vertices maintain PST?}\label{sec:wvmPST} 
We keep the same vertex ordering as before, where vertices of the 4-cube are labeled as in Fig.~1 and for the $n$-cube in accordance with the dictionary ordering of the vertices of the Cartesian product. Let  $S=\{1+2^4m_4+\cdots+2^{n-1}m_{n-1}, 16+2^4m_4+\cdots+2^{n-1}m_{n-1}\,|\,  m_k\in\{0,1\} \text{ for } k=4,\dots,n-1\}$.

\vspace*{12pt}
\noindent
{\bf Theorem~4:}\label{paritialPST} Let $n>4$. Then for any given partially switched $n$-cube,
%There exists PST in any partially switched $n$-cubes, and, for any $n>4$, all of the partially switched $n$-cubes have the same PST vertex pairs. 
%Furthermore, 
at least 1/8 of its vertices pair up to exhibit PST at time $\pi/2$. 
Specifically, for any partially switched $n$-cube, all the vertices in the set $S$  pair up to exhibit PST, with PST   vertex pairs $1+2^4m_4+2^5m_5+\cdots+2^{n-2}m_{n-2}+2^{n-1}m_{n-1}$ and 
$16+2^4(1-m_4)+2^5(1-m_5)+\cdots+2^{n-2}(1-m_{n-2})+2^{n-1}(1-m_{n-1})$,
where $m_j\in\{0,1\}$  for $j=4,\ldots,n-1$. 
%for any vertex $v$ in the set $S=\{1+2^4m_4+\cdots+m_{n-1}\times (2^{n-1}), 16+2^4m_4+\cdots+m_{n-1}\times (2^{n-1})\,|\,  m_k\in\{0,1\} \text{ for } k=m_4,\cdots,m_{n-1}\}$ (the set of singleton vertices in the partition $\pi_1$), there is vertex $u\in S$, such that the fidelity of state transfer between $u$ and $v$ is 1 for any partially switched $n$-cube.
Furthermore, at any time $t$, the fidelity of state transfer from vertex $j$ to any other vertex is the same for any partially switched $n$-cube as it is for the $n$-cube. 
\vspace*{12pt}

% however, the proof below relies   on considering the corresponding entries in the unitary matrix $e^{i(\pi/2)A}$, avoiding the use of quotient graphs. 
Theorem 2 in \cite{quogra}
states that if we have an equitable partition with $u, v$ as singleton cells, then the fidelity from $u$ to $v$ at any time $t$ is the same in the original graph as it is in the symmetrized quotient graph. Therefore one could use the fact that the $n$-cube and a partially switched $n$-cube  have the same quotient graph according to the partition $\pi_1$ as mentioned above Theorem~\ref{thm:switchedPSTpairs} to prove the PST pairs results.
But our arguments show something stronger, namely  that for the $n$-cube, switched $n$-cube and partially switched $n$-cubes (and, as we show in Propositions 10 and 13, convex combinations and time-switched systems), at any time $t$, the fidelity from vertex 1 to every other vertex (not just the fidelity to vertex $2^n$) is preserved. So, there is more information than what we can deduce  from using Theorem 2 in \cite{quogra}.

\noindent
{\bf Proof:}
%First we show that for any positive integer $k$, $A_{4,1}^k$ and $A_{4,2}^k$ have the same first row and last row (16-th row).
%We prove it for the first row, and the result for the 16-th row follows in exactly the same way.
By direct computation, we know $\bra{1}A_{4,1}^k=\bra{1}A_{4,2}^k$ for $k=1,2,3,4$. 
Since $A_{4,1}$ and $A_{4,2}$ have the same minimal polynomial $x^5-20x^3+64x$,
we know $\bra{1}A_{4,1}^k=\bra{1}A_{4,2}^k$ for any positive integer $k$, i.e., for any positive integer $k$, $A_{4,1}^k$ and $A_{4,2}^k$ have the same first row. Similarly, they have the same $16$-th row.
Therefore, for any nonnegative integers 
$j_1,j_2,\dots, j_{2s}$, the matrix $A_{4,1}^{j_1}A_{4,2}^{j_2}A_{4,1}^{j_3}\cdots A_{4,2}^{j_{2s}}$ has the same first row 
as $A_{4,1}^{j_1+\cdots +j_{2s}}$, since
\begin{eqnarray*}
\bra{1}A_{4,1}^{j_1}A_{4,2}^{j_2}A_{4,1}^{j_3}\cdots A_{4,2}^{j_{2s}}
&=&\bra{1}A_{4,2}^{j_1}A_{4,2}^{j_2}A_{4,1}^{j_3}\cdots A_{4,2}^{j_{2s}}=\bra{1}A_{4,2}^{j_1+j_2}A_{4,1}^{j_3}\cdots A_{4,2}^{j_{2s}}\\
&=&\bra{1}A_{4,1}^{j_1+j_2+j_3}\cdots A_{4,2}^{j_{2s}}
=\cdots%= \bra{1}A_{4*}^{j_1+j_2+j_3+j_\ell}
\\&=&
\bra{1}A_{4,1}^{j_1+j_2+j_3+\cdots+j_{2s}}
\\&=&
\bra{1}A_{4,2}^{j_1+j_2+j_3+\cdots+j_{2s}}.
\end{eqnarray*}
%Similarly, $\bra{1}A_{4,1}^{j_1}A_{4,2}^{j_2}A_{4,1}^{j_3}\cdots A_{4,2}^{j_{2s}}=\bra{1}A_{4,2}^{j_1+j_2+j_3+\cdots+j_{2s}}$. 
As mentioned earlier, for $n> 4$, the adjacency matrix of a partially switched $n$-cube is of the form
$A_{n,p}=\diag(A_{4*},\cdots,A_{4*})+A_{n-4,1}\otimes I_4$, where $A_{n-4,1}$ is the   adjacency matrix of the $(n-4)$-cube, and $*$ represents 1 or 2.
For any positive  integer $k$, 
%For the unitary fidelity matrix $Up=e^{i\pi/2 Ap}=\sum_{j=0}^{\infty}\frac{(i\pi/2 Ap)^j}{j!}$, 
each block of the matrix $A_{n,p}^k$ is of the form $\sum c(j_1,j_2,\ldots,j_{2s})A_{4,1}^{j_1}$ $ A_{4,2}^{j_2}\cdots A_{4,2}^{j_{2s}}$ for some nonnegative integers
$j_1,\dots, j_{2s}$ and some real number $c(j_1,j_2,\ldots,j_{2s})$, which has the same first row as $\sum c(j_1,j_2,\ldots,j_{2s})A_{4,1}^{j_1+\cdots+j_{2s}}$, the corresponding block in $A_{n,1}^k$.
%$U=e^{i\pi/2 A_{n,1}}$. 
%Therefore $\bra{j}Up=\bra{j}U$ for $j=1,16,2^4+1,2^4+16,\cdots$.
Therefore  $\bra{\ell}A_{n,p}^k=\bra{\ell}A_{n,1}^k$ for $\ell=1,16,2^4+1,2^4+16, 2^5+1,2^5+16, 2^5+2^4+1, 2^5+2^4+16, \dots$ (these rows correspond to the first and  16-th vertices in each copy of the 4-cube or switched 4-cube).
Hence for the unitary matrices $U_{n,p}(t)=e^{it A_{n,p}}=\sum_{j=0}^{\infty}\frac{(it A_{n,p})^j}{j!}$ for $A_{n,p}$ and $U_{n,1}(t)=e^{it A_{n,1}}$ for $A_{n,1}$, 
$\bra{\ell}U_{n,p}(t)=\bra{\ell}U_{n,1}(t)$ for these $\ell$ and any time $t$.
As a result, at any time $t$, each of the vertices in the set $S$ has the same probability of state transfer to any other vertex as it has in the $n$-cube.
In the $n$-cube, there is PST between any two vertices of distance $n$ at time $t=\pi/2$, which correspond to 
vertices $k$ and $2^n+1-k$ i.e., $|\bra{k}U_{n,1}(\pi/2)\ket{2^n+1-k}|=1$. Therefore, in any partially switched $n$-cube,
there is PST between vertices $1+2^4m_4+2^5m_5+\cdots+2^{n-2}m_{n-2}+2^{n-1}m_{n-1}$ and 
$16+2^4(1-m_4)+2^5(1-m_5)+\cdots+2^{n-2}(1-m_{n-2})+2^{n-1}(1-m_{n-1})$ at time $\pi/2$,
% between vertices $16+(2^4)^{u_4}+(2^5)^{u_5}+\cdots+(2^{n-2})^{u_{n-2}}$ and $1+(2^4)^{w_4}+(2^5)^{w_5}+\cdots+(2^{n-2})^{w_{n-2}}+2^{n-1}$,
where $m_j\in\{0,1\}$  for $j=4,\ldots,n-1$. %We note that it is still true that the sum of the indices for PST pairs is $2^n+1$, as it is for the $n$-cube.
%{\color{red} haven't figured out how to prove there are no other PST pairs}
$\square$

 Below we conjecture that the lower bound of 1/8 of the vertices in Theorem~4 is in fact exact (that is, exactly 1/8 of the vertices of partially switched $n$-cubes pair up to exhibit PST at time $\pi/2$). As a motivating example, we consider $A_{5,3}$, the first interesting partially switched $n$-cube, and verify that the bound is attained in this case. The verification process is rather tedious, but we include the technical details for completeness. The main takeaway is that the smallest nontrivial example does indeed attain our  lower bound.

We first recall some spectral graph theory results. The adjacency matrix $A(G)$ of a graph $G$ is real and symmetric, so it has a spectral decomposition $A(G)=\sum_{r=1}^s \lambda_rE_r$, where  $\lambda_1,\ldots, \lambda_s$ are all the distinct eigenvalues of $A(G)$, and $E_j$  represents the orthogonal projection onto the eigenspace associated with eigenvalue $\lambda_j$.  Given a vertex $u\in V(G)$, its characteristic (indicator) vector is $\ket{u}\in\mathbb{R}^{|V(G)|}$. The eigenvalue support of $\ket{u}$ is defined to be the set of eigenvalues $\lambda_r$ of $A(G)$, such that $E_r\ket{u}\neq0$.

\vspace*{12pt}
\noindent
{\bf Theorem~5:}\cite[Theorem 3.1]{whenPST}\label{Th:eigsup}
Let $X$ be a graph and $u$ be a vertex in $X$ at which $X$ is periodic. If $\theta_k$, $\theta_{\ell}$, $\theta_r$, $\theta_s$ are eigenvalues in the support of $\ket{u}$ and 
$\theta_r\neq\theta_s$, then $\frac{\theta_k-\theta_{\ell}}{\theta_r-\theta_s}\in \mathbb{Q}$. (Therefore if there are two integer eigenvalues in the support of $u$, then all the eigenvalues in the support of $u$ are integers.)
%Suppose $G$ is a connected graph with at least three vertices where perfect state transfer from $u$ to $v$ occurs at time $t=t_0$. Let $\delta$ be the dimension of the $A(G)$-invariant subspace generated by $\ket{u}$. Then $\delta\geq3$ and either all the eigenvalues in the eigenvalue support of $u$ are integers, or they are all of the form $\frac12(a+b_j\sqrt{\Delta})$ where $\Delta$ is a square-free integer and $a$ and $b_j$ are integers.
\vspace*{12pt}

\vspace*{12pt}
\noindent
{\bf Remark~6:}  Let $G$ be a graph on $m$ vertices, and  $u$ be a vertex of $G$. Then the eigenvalue $\lambda_r$ of $A(G)$ is in the eigenvalue support of $\ket{u}$ if we can find a (normalized) eigenvector $\ket{v_1}$ of $A(G)$ associated to $\lambda_r$, such that $\braket{v_1}{u}\neq 0$. %Since for any eigenvector $v_1$ of $A(G)$ associated to $\lambda_r$, we can extend it to
From $\ket{v_1}$, we can get   a basis $\{\ket{v_1}, \ket{v_2},\ldots,\ket{v_k}\}$ of the eigenspace associated to $\lambda_r$, then by the Gram-Schmidt procedure,  we can get an orthonormal basis $\ket{w_1},\ldots, \ket{w_k}$ of the eigenspace. Now we have $E_r\ket{u}=(\ketbra{w_1}{w_1}+\ketbra{w_2}{w_2}+\cdots+\ketbra{w_k}{w_k})\ket{u}=\braket{w_1}{u}\ket{w_1}+\cdots+\braket{w_k}{u}\ket{w_k}$ for any vertex $u$,  and $E_r\ket{u}=0$ if and only if $\braket{w_j}{u}=0$ for all $j=1,\ldots,k$, i.e., all the eigenvectors of $A(G)$ associated to $\lambda_r$ have their $u$-th entry  equal to 0. This implies that for any eigenvalue $\lambda_r$ of $A(G)$, if it has a corresponding eigenvector
 whose $u$-th entry is not 0, then $\lambda_r$ is in the eigenvalue support of $\ket{u}$.
 \vspace*{12pt}

We are now in the position to consider the example of $A_{5,3}$. 

\vspace*{12pt}
\noindent
{\bf Example~7:}\label{A53PST}
{\rm Consider $A_{5,3}=\begin{bmatrix}A_{4,1}&I_4\\I_4&A_{4,2}\end{bmatrix}$,
it has $\lambda_1=5$ as a simple eigenvalue, with $v_1=\mathbf{1}_5$ being a corresponding eigenvector. %Therefore $E_1\ket{u}=\frac{1}{32}\mathbf{1}_{(32)}\mathbf{1}_{(32)}^T\ket{u}=\frac{1}{32}\mathbf{1}_{(32)}\neq0$ for any $u=1,\ldots, 32$, which means 
%Since all the entries of $\mathbf{1}$ are nonzero, we know $\lambda_1=5$ is in the eigenvalue support of all the vertices of $G(A_{5,3})$. 
By direct computation, we know $v_2=\begin{bmatrix}\mathbf{1}_4\\-\mathbf{1}_4\end{bmatrix}$ is an eigenvector of $A_{5,3}$ associated to eigenvalue 3. Since all the entries of $v_1$ and $v_2$ are nonzero, from the above remark, we know $\lambda_1=5$ and $\lambda_2=3$ are both in the eigenvalue support of all the vertices.
\iffalse
For the eigenvalue $\lambda_3=2.7093$ of $A_{5,3}$, %its associated eigenprojection matrix is $E_3=v_1v_1^T+v_2v_2^T+v_3_3^T$, where 
$v_1=[0.0000,  0.1195,   {-0.1379},   0.0651, {-0.0467},   {-0.0342},   0.3424,   {-0.1350},    0.1350,$ $  {-0.3424},    0.0342,   0.1014,  {-0.1414},    0.2993, 
{-0.2593},  {-0.0000},    0.0000,   {-0.1195},  $ $ 0.1379,   {-0.0651},    0.0467,   {-0.0342},    0.3424,   {-0.1350},    0.1350,   {-0.3424},    0.0342,$ $    0.1014,
 {-0.1414},    0.2993,   {-0.2593},   {-0.0000}]^T$ is an associated eigenvector, using the above remark we can see that $\lambda_3=2.7093$ is in the eigenvalue support of every vertex $j\in Z=\{1,2,\ldots,32\}\backslash \{1,16,17,32\}$. From Theorem~5 we know that no vertices in the set $Z$ exhibits PST. Combing this result with Theorem~4, we know set of vertices of $A_{5,3}$ with PST are exactly $\{1,16,17,32\}$.

%Every root of the polynomial 
Assume $\lambda_1\geq\cdots\geq\lambda_6$ are the roots of $p(\lambda)=-\lambda^6+11\lambda^4-27\lambda^2+1=-(\lambda^3+\lambda^2-5\lambda-1)(\lambda^3-\lambda^2-5\lambda+1)$; they are all eigenvalues  of $A_{5,3}$. Each of them is irrational (non-integer eigenvalues of any monic integer-coefficients polynomial  are irrational), with minimal polynomial  either $(\lambda^3+\lambda^2-5\lambda-1)$ or $(\lambda^3-\lambda^2-5\lambda+1)$. For example, the approximate value of the largest root is $\lambda_1=2.7093$.
\fi 

Now let $p(\lambda)=-\lambda^6+11\lambda^4-27\lambda^2+1=-(\lambda^3+\lambda^2-5\lambda-1)(\lambda^3-\lambda^2-5\lambda+1)$, then it has 6 real roots, for example, by the Intermediate Value Theorem we can check it has a root $\lambda_3$ between 2.7 and 2.8. Assume the 6 roots are $\lambda_3\geq\cdots\geq\lambda_8$; they are all eigenvalues of $A_{5,3}$, and each of them is irrational (non-integer roots of a monic integer-coefficients polynomial are irrational), with minimal polynomial either $(\lambda^3+\lambda^2-5\lambda-1)$ or $(\lambda^3-\lambda^2-5\lambda+1)$. For $k=3,\ldots,8$, the eigenvalue $\lambda_k$ has an associated 
%The eigenvalue $\lambda=\lambda_k$, $k=1,\ldots,6$, has an associated 
eigenvector 
\iffalse %%%%%%%%
$v(\lambda_k)= [0,            2\lambda_k(\lambda_k^4 - 10\lambda_k^2 + 17),       -6\lambda_k(\lambda_k^4 - 10\lambda_k^2 + 17),
  2\lambda_k(\lambda_k^4 - 10\lambda_k^2 + 17),          2\lambda_k(\lambda_k^4 - 10\lambda_k^2 + 17),
8\lambda_k^2 - 8,           8 - 8\lambda_k^2,           8\lambda_k^2 - 8,      8 - 8\lambda_k^2,        8\lambda_k^2 - 8,       8 - 8\lambda_k^2,
2\lambda_k(\lambda_k^2 - 5)^2,     2\lambda_k(\lambda_k^2 - 5)^2,        -6\lambda_k(\lambda_k^2 - 5)^2,
 2\lambda_k(\lambda_k^2 - 5)^2,      0,       0,       2\lambda_k^4 - 12\lambda_k^2 - 6,        36\lambda_k^2 - 6\lambda_k^4 + 18,
2\lambda_k^4 - 12\lambda_k^2 - 6,       2\lambda_k^4 - 12\lambda_k^2 - 6,          8\lambda_k(\lambda_k^2 - 5),
-8\lambda_k(\lambda_k^2 - 5),         8\lambda_k(\lambda_k^2 - 5),      -8\lambda_k(\lambda_k^2 - 5),
8\lambda_k(\lambda_k^2 - 5),        -8\lambda_k(\lambda_k^2 - 5),       2\lambda_k^4 - 12\lambda_k^2 + 10,
2\lambda_k^4 - 12\lambda_k^2 + 10,       36\lambda_k^2 - 6\lambda_k^4 - 30,      2\lambda_k^4 - 12\lambda_k^2 + 10,     0]^T$. 
\fi%%%%%%%%%%%
$v(\lambda_k)=[0,a,-3a,a,a,b,-b,b,-b,b,-b,c,c,-3c,c,0,0,d,-3d,\\
d,d,e,-e,e,-e,e,-e,d+4,d+4,-3(d+4),d+4,0]^T$,
where $a=2\lambda_k(\lambda_k^4 - 10\lambda_k^2 + 17)$, $b=8\lambda_k^2 - 8$, $c=2\lambda_k(\lambda_k^2 - 5)^2$, $d=2\lambda_k^4 - 12\lambda_k^2 - 6$,  and   $e=8\lambda_k(\lambda_k^2 - 5)$.
Note that for each $u\in Z=\{1,2,\ldots,32\}\backslash \{1,16,17,32\}$, the entry $v(\lambda_k)_u$ is not divisible by the minimal polynomial of $\lambda_k$, and therefore none of these entries are zero. Again from the above  remark, we know that for each $k=3,\ldots,8$, $\lambda_k$  is in the eigenvalue support of every vertex $u\in Z$. 
Now for each $u\in Z$, $\lambda_1=5$, $\lambda_2=3$, and $\lambda_3\in[2.7,\; 2.8]$ are in the eigenvalue support of $\ket{u}$.
%But there are no integers $a, b_1, b_2, b_3$ and a square-free integer $\Delta$ such that $\lambda_{\ell}=\frac12(a+b_{\ell}\sqrt{\Delta})$, for $\ell=1,2,3$.
Since PST implies periodicity, Theorem~5 implies that, no vertices in the set $Z$ exhibit PST.  Combining this result with Theorem~4, we know the set of vertices of $A_{5,3}$ with PST is exactly $\{1,16,17,32\}$.}
\vspace*{12pt}

\vspace*{12pt}
\noindent
{\bf Conjecture~8:} Let $n>4$. For any partially switched $n$-cube, the set of its vertices that exhibits PST is exactly the set $S$ we give in Theorem~4, therefore exactly 1/8 of the vertices of a partially switched $n$-cube pair up to have PST.
\vspace*{12pt}

\subsection{Other Variants}
A dual-channel quantum  directional  coupler was introduced in \cite{routingNik} as a means to selectively transfer a state to either of the two output ports in a controlled
and deterministic way; the Hamiltonian describing these dynamics is $\mathcal H=\mathcal H_h+\mathcal H_v$ where $\mathcal H_h$ describes the energy related to the source or drain channel, while $\mathcal H_v$ describes the energy of the inter-channel dynamics. Quantum state transfer between vertices in parallel, multi-user networks, is described in \cite{routingCS}, where one sender-receiver pair uses each channel at a time to achieve optimal routing. A description of how to design a large quantum network out of smaller independent subsystems is described in \cite{PRK}. Motivated by this literature on quantum routing as a technique in quantum state transfer, we consider several ways to oscillate between the graphs considered herein while preserving PST for many vertex pairs, and in fact maintaining all the dynamics of fidelity (namely, the fidelity function is exactly the same) for many vertex pairs.

We first consider different convex combinations of each 4-cube or switched 4-cube block of a partially switched $n$-cube. General linear combinations could also be used, though the PST time would change. 

\vspace*{12pt}
\noindent
{\bf Remark~9:}\label{Rem:weig} {\rm We can generalize partially switched $n$-cubes to specially weighted ones. Consider a convex combination  of the 4-cube and the switched 4-cube. The resulting graph $G$ has adjacency matrix $M=pA_{4,1}+(1-p)A_{4,2}$ for $0\leq p\leq 1$. Using the same techniques we used in the proof of Theorem~4, we can see there is  perfect state transfer between vertex 1 and vertex 16. Furthermore, by induction we can see that for the graph %$Gm=Q_{n-4}\square G$ has PST pairs as stated in Theorem  ~4. The 
$F=Q_{n-4}\square G$, whose adjacency matrix %of $Gm$ is 
is $A(F)=I_{n-4}\otimes M+C_{n-4}\otimes I_4=\diag(M,\cdots, M)+C_{n-4}\otimes I_4$, every vertex in the set $S$ as mentioned in  Theorem~4 exhibits PST. %has the same PST pairs as those stated in Theorem  ~4. %Now we can generalize the graph again. Consider the 
A similar statement holds for the graph $\tilde{F}$ with (nonnegative) adjacency matrix $A=\diag(M_1,M_2,\cdots, M_{2^{n-4}})+C_{n-4}\otimes I_4$, where $M_j=p_jA_{4,1}+(1-p_j)A_{4,2}$ for $0\leq p_j\leq 1, j=1,\ldots, 2^{n-4}$. % As above,  we can show that the corresponding weighted graph has the same PST pairs as those stated in Theorem~4. 
Note that this new family of graphs contains all the other cubes as special cases: when $p_1=\cdots=p_{2^{n-4}}=1$, we have the $n$-cube, where all the vertices pair up to exhibit perfect state transfer; when $p_1=\cdots=p_{2^{n-4}}=0$, we have the switched $n$-cube, where exactly half of the vertices pair up to exhibit perfect state transfer; when $p_1,\ldots, p_{n-4}\in\{0,1\}$ and not all of them are equal, then we have a partially switched $n$-cube, and Theorem~4 gives a list of vertex pairs having PST.

\vspace*{12pt}
More generally, we can consider convex combinations of arbitrary graphs on $m$ vertices whose adjacency matrices satisfy some specific conditions for some row. 
%where for any nonnegative integer $k$, the $k$-th power of their adjacency matrices have the same $j$-th row for some $1\leq j\leq m$. Without lose of generality, we assume $j=1$.
%$pA(G_1)+(1-p)A(G_2)$ for  $0\leq p\leq 1$.

\vspace*{12pt}
\noindent
{\bf Proposition~10:}\label{Prop:convex}
Let $G_1,G_2,\ldots, G_k$ be graphs on $m$ vertices, whose corresponding adjacency matrices are $A(G_1),A(G_2),\ldots,A(G_k)$, respectively. 
Suppose that for some $\ell\in \{1, \dots, m\}$, $\bra{\ell}A(G_{r})^j=\bra{\ell}A(G_{s})^j$ for every positive integer $j$ and any $r,s=1,\ldots,k$. If there is PST in any one of the $k$ graphs from vertex $\ell$ to some other vertex $u$ at time $t=t_0$, then all the other graphs have PST between vertex $\ell$ and $u$ at time $t_0$, as well as the weighted graph $G$ with adjacency matrix $A=c_1A(G_1)+\cdots+c_kA(G_k)$, where  $0\leq c_r\leq 1$ for $r=1,\ldots,k$,  and $c_1+\cdots +c_k=1$.

\vspace*{3pt}
% is equal to the first row of  $e^{itA(G_2)}$ for every time $t\in \R^+$. Then the graph $G_3$ whose adjacency matrix is the convex combination $pA(G_1)+(1-p)A(G_2)$, 
The argument is similar to that given in Remark 9 as well as the proof of Theorem~4.  
In particular, we note that $A^j=(c_1A(G_1)+\cdots+c_kA(G_k))^j$ has the same $\ell$-th row as $A(G_r)^j=(c_1A(G_r)+\cdots+c_kA(G_r))^j$ for any nonnegative integer $j$ and $r=1,\ldots,k$.

\vspace*{12pt}
\noindent
{\bf Corollary~11:}\label{cor:convex}
Any convex combination of the $n$-cube, the switched $n$-cube, some partially switched $n$-cube, and the weighted matrices in Remark 9,
%$pA_{n,1}+(1-p)A_{n,2}$ for $0\leq p\leq 1$ will 
has PST between vertices  $1+2^4m_4+2^5m_5+\cdots+2^{n-2}m_{n-2}+2^{n-1}m_{n-1}$ and 
$16+2^4(1-m_4)+2^5(1-m_5)+\cdots+2^{n-2}(1-m_{n-2})+2^{n-1}(1-m_{n-1})$, %$1+2^4m_4+2^5m_5+\cdots+m_{n-2}\times(2^{n-2})+m_{n-1}\times(2^{n-1})$ and 
%$16+(1-m_4)\times(2^4)+(1-m_5)\times(2^5)+\cdots+(1-m_{n-2})\times(2^{n-2})+(1-m_{n-1})\times(2^{n-1})$,
 at time $\pi/2$, where $m_j\in\{0,1\}$  for $j=4,\ldots,n-1$. %Moreover, any partially switched $n$-cube whose $n$-cube and switched $n$-cube blocks are replaced with convex combinations of the $n$-cube and switched $n$-cube will also have PST between vertices??? at time $\pi/2$. 
\vspace*{12pt}

\vspace*{12pt}
\noindent
{\bf Remark~12:}  We already know that any convex combination of $A_{4,1}$ and $A_{4,2}$ exhibits PST between vertex 1 and 16. Here we give some spectrum properties of such convex combinations.

For any $0\leq p\leq 1$, the weighed graphs with adjacency matrices $pA_{4,1}+(1-p)A_{4,2}$ and $pA_{4,2}+(1-p)A_{4,1}$, respectively, are cospectral to each other, with the similarity matrix being the symmetric orthogonal matrix $Q$  given by $QA_{4,1}Q=A_{4,2}$ and $QA_{4,2}Q=A_{4,1}$ ($Q$ exists by  \cite{GM}). Alternatively, let $C=1/2A_{4,1}+1/2A_{4,2}$ and $E=A_{4,2}-A_{4,1}$. Then for any $0<\alpha\leq1/2$, the two nonnegative matrices $C+\alpha E$ and $C-\alpha E$ have the same spectrum (indeed, since $QCQ=C$ and  $QEQ=-E$, we have $Q(C+\alpha E)Q=C-\alpha E$).
The eigenvalues of $C+\alpha E$ and $C-\alpha E$  are $\pm4$ (with multiplicity 1), $\pm 2$ (with multiplicity 1), 0 (with multiplicity 6), and $\sqrt{2+8\alpha^2}$ (with multiplicity 3), which can be checked by calculating the rank of the corresponding matrices.

Similarly, for the adjacency matrix $C_n=I_{n-4}\otimes A_{4,1}+C_{n-4}\otimes I_4$ of $Q_n$, if we replace the diagonal blocks $A_{4,1}$ by different convex combinations of $A_{4,1}$ and $A_{4,2}$, %we have for $0\leq p_j\leq 1, j=1,\ldots, n-4$, 
then the nonnegative matrices $\diag(p_1A_{4,1}+(1-p_1)A_{4,2},p_2A_{4,1}+(1-p_2)A_{4,2},\cdots, p_{2^{n-4}}A_{4,1}+(1-p_{2^{n-4}})A_{4,2})+C_{n-4}\otimes I_4$, and   $\diag(p_1A_{4,2}+(1-p_1)A_{4,1},p_2A_{4,2}+(1-p_2)A_{4,1},\cdots, p_{2^{n-4}}A_{4,2}+(1-p_{2^{n-4}})A_{4,1})+C_{n-4}\otimes I_4$ have the same spectrum (similar through the matrix $\diag(Q,Q,\ldots, Q)$), where $0\leq p_j\leq 1, j=1,\ldots, 2^{n-4}$.

\vspace*{12pt}

As another variant, we consider switched systems where one employs a switching function to change between systems at particular times (this can be done in the absence of GM switching---it is a coincidence in naming). For example, one might use the spin network associated to the hypercube from time $t=0$ to time $t=t_1$, then change to the spin network associated to the switched cube from time $t=t_1$ to time $t=t_2$, change to use a partially switched hypercube from time $t=t_2$ to time $t=t_3$, and so on, up to time $t_r=\pi/2$, when the $n$-cube, the switched $n$-cube, and any partially switched $n$-cube have PST. We show that this new system (whose Hamiltonian changes with respect to time) has PST for vertices in the set $S$ as mentioned in Theorem~4.
%common PST vertex pairs (that is, employing a switch system preserves PST for the PST vertex pairs that the graphs have in common). 

The motivation here is potential stability issues in the lab: spin networks are created in the lab with magnets and other devices and may be unstable, especially for long periods of time. Thus, one might wish to send a state along the first network until one loses confidence in the stability, then one can change to the second network and continue sending the state through this ``fresh'' network while rebooting the first. This would be an example of a quantum state transfer protocol requiring external modulation; such external modulation approaches typically increase the effectiveness of the state transfer, but it may be undesirable to use a protocol that relies heavily on a ``hands on'' approach. A binary switching between spin networks may be a useful compromise. Our approach is motivated by switched systems in control theory; see, e.g.\ \cite{Dan}.

\vspace*{12pt}
\noindent
{\bf Proposition~13:}\label{Prop:switching}
Assume $r$ is some positive integer. For $j=1,\ldots,r$, let $G_j$ % and $G_2$ 
be any of the following: the $n$-cube, the switched $n$-cube, the partially switched $n$-cubes, or convex combinations as described in Corollary~11. If a quantum state is transferred through the network $G_1$ for $0\leq t\leq t_1$, $G_2$ for $t_1\leq t\leq t_2$, $G_3$ for $t_2\leq t\leq t_3$, $\dots$, $G_r$  for $t_{r-1}\leq t\leq t_r=\pi/2$, %with $t_r=\pi/2$, (and $i=1$ for $r$ odd and $i=2$ for $r$ even) 
then the quantum system with time-dependent Hamiltonian $\mathcal H_t$ 
is guaranteed to have PST at time $\pi/2$ for the vertices in the set $S$ as mentioned in Theorem~4.
%\vspace*{12pt}
Further,  if for each $j=1,\ldots, r$, $G_j$ is either the $n$-cube or the switched $n$-cube,  and in addition  at least one $G_j$ is the switched $n$-cube, then the set of vertices exhibiting PST in this system is exactly the set of vertices exhibiting PST in the switched $n$-cube (one half of all the vertices).

\noindent
{\bf Proof:} 
We consider the case $r=2$, the general case follows from induction. Fix a vertex $\ell\in S$.
We have
\begin{eqnarray*}
%\bra{1}\exp\{it_1A_{n,1}+it_2A_{n,2}\}&=&
&&\bra{\ell}\exp\{it_1A(G_1)\}\exp\{i(t-t_1)A(G_2)\}\\
&=&\bra{\ell}\sum_{k=0}^\infty \frac{(it_1)^k A(G_1)^k}{k!}\exp\{i(t-t_1)A(G_2)\}\\
&=&\bra{\ell}\sum_{k=0}^\infty \frac{(it_1)^k A(G_2)^k}{k!}\exp\{i(t-t_1)A(G_2)\}\textnormal{ by Proposition~10}\\
&=&\bra{\ell}\exp\{it_1A(G_2)\}\exp\{i(t-t_1)A(G_2)\}\\
&=&\bra{\ell}\exp\{i(t_1+t-t_1)A(G_2)\}=\bra{\ell}\exp\{itA(G_2)\}.
\end{eqnarray*}
Thus the problem reduces to finding PST pairs in $S$ for $A(G_2)$. 
$\square$

We will analyse the variants discussed above in terms of their   sensitivity to readout time errors in Section \ref{sec:analysis}. 
When there is PST these variants have the same sensitivity to readout time errors as the original hypercube when considering PST pairs from the set $S$ defined above Theorem~4.  

\section{Sensitivity analysis with respect to readout time errors}\label{sec:analysis}
Recall $S=\{1+2^4m_4+\cdots+2^{n-1}m_{n-1}, 16+2^4m_4+\cdots+2^{n-1}m_{n-1}\,|\,  m_k\in\{0,1\} \text{ for } k=4,\dots,n-1\}$. The sensitivity of the probability (fidelity) of state transfer  with respect to readout time is typically analyzed through the first derivative. An analysis of the $k$th derivatives (for any $k\in \mathbb{N}$) for weighted graphs with PST was done in \cite{Steve2015}. Here, we consider both the first and second derivatives.

\vspace*{12pt}
\noindent
{\bf Theorem~14:} The $n$-cube, the switched $n$-cube, the partially switched $n$-cubes, and the other $n$-cube variants discussed herein all have the same derivatives with respect to time $t$ at time $t=\pi/2$ %sensitivity with respect to readout time errors 
for the PST pairs of vertices in the set $S$. 

\vspace*{12pt}
\noindent
{\bf Proof:} 

As in Theorem~4, we have already shown that for any vertex $j\in S$, there is PST between vertex $j$ and vertex $2^n+1-j$ at time $\pi/2$ for the (switched, partially switched) $n$-cube, and %if there is PST from vertex $j$ to $k$ in a partially switched $n$-cube, then 
$\bra{j}U_{n,p}(t)=\bra{j}U_{n,1}(t)=\bra{j}U_{n,2}(t)$. Therefore the fidelity of state transfer from vertex $j$ to any other vertex $k$ is the same as it is 
in the $n$-cube at any time $t$. It follows that the three types of $n$-cubes have the same derivatives with respect to readout time $t$ %sensitivity to timing errors 
at time $t=\pi/2$. For an undirected graph $G$ exhibiting PST between vertices $r$ and $s$ at time $t=t_0$, 
the derivatives of fidelity with respect with readout time $t$ at time $t=t_0$ is given in \cite{Steve2015}:
\begin{eqnarray*}
\frac{d^kp}{dt^k}\Bigg|_{t=t_0}=\begin{cases}
    (-1)^{(k\mod 4)/2}  \sum_{\ell=0}^k  (-1)^\ell{k\choose \ell}\bra{s}\mathcal{H}^\ell\ket{s}\bra{s}\mathcal{H}^{k-\ell}\ket{s}   & \quad \text{if } k \text{ is even}\\
    0 & \quad \text{if } k \text{ is odd.}
  \end{cases}
\end{eqnarray*} 
From this we find that, for the PST vertex pairs in the set $S$,  $\frac{dp}{dt}|_{t=\frac{\pi}{2}}=0$, and $\frac{d^2p}{dt^2}|_{t=\frac{\pi}{2}}=-2n$, where $p$ is the fidelity of state transfer at time $t$ between PST vertex pairs in $S$.
Similarly, we can use the proof in Proposition~10 to prove this result for convex combinations, and use Proposition~13 to prove it for the switched system. 
$\square$

\section{Conclusion} 
The hypercube and the more general class of cubelike graphs whose elements in the connection set have nonzero sums have seen much attention recently as they have been shown to exhibit PST between pairs of vertices at time $\pi/2$, where the PST pairing is determined by the sum of the elements in the connection set. Here, we perform various perturbations on the hypercube while maintaining PST for a subset of vertices, including a perturbation that allows for a time-dependent Hamiltonian, which may be of practical use. The fidelity of state transfer involving vertices where PST occurs in our various perturbed hypercubes have the same sensitivity to readout time errors as the original hypercube, thus identifying infinite families of graphs sharing the highly desirable properties of PST and maximal PST distance. It would be of interest to see if further perturbations can be done, or how radical the perturbations can be, before completely losing the property of PST.

\section*{Acknowledgements}
\noindent
S.K.~and S.P.~were supported by NSERC Discovery Grants RGPIN/6123-2014  and 1174582, respectively;  S.P.~is also supported by the Canada Foundation for Innovation (CFI) grant number 35711, and the Canada Research Chairs (CRC) Program grant number 231250.
  X.Z.\ was supported by the University of Manitoba's Faculty of Science and Faculty of Graduate Studies.


\begin{thebibliography}{000}


\bibitem{Cubes} M.\ Christandl, N.\ Datta, T.\ Dorlas, A.\ Ekert, A.\ Kay, and A.\ Landahl, \emph{Perfect transfer of arbitrary states in quantum spin networks}, Phys.\ Rev.\ A \textbf{71} (2005), 032312.
\bibitem{Cubelike} W.-C.\ Cheung, C.\ Godsil, \emph{Perfect state transfer in cubelike graphs}, Linear Algebra Appl.\ \textbf{435} (2011), 2468-2474.
\bibitem{MBHadamards} N.~Johnston, S.\ Kirkland, S.\ Plosker, R.\ Storey, and X.\ Zhang, \emph{Perfect quantum state transfer using Hadamard diagonalizable weighted graphs}, Linear Algebra and its Applications, \textbf{531} (2017), 375-398.
\bibitem{BFK} S.\ Barik, S.\ Fallat, and S.\ Kirkland, \emph{On Hadamard diagonalizable graphs}. Linear Algebra Appl.\ \textbf{435} (2011), 1885-1902.

\bibitem{potentials} A.~Casaccino, S.~Lloyd, S.~Mancini,  and S.~Severini, \emph{Quantum state transfer through a qubit network with energy shifts and fluctuations},  Int.~J.\ Quantum Inf.\ \textbf{7}(08) (2009), pp.\ 1417-1427.
\bibitem{KemptonPST} M.~Kempton, G.~Lippner, and S.-T.~Yau, \emph{Perfect state transfer on graphs with a potential}, Quantum Inf.\ Comput.\ \textbf{17} no.\ 3\& 4 (2017), pp.\ 303--327.
\bibitem{RV} S.~Kirkland,D.~McLaren, R.~Pereira, S.~Plosker, and X.~Zhang, \emph{Perfect quantum state transfer in weighted paths with potentials (loops) using orthogonal polynomials},  Linear Multilinear A.\ (2018), 1-19.
\bibitem{CDDEKL} M.~Christandl, N.~Datta,  T.~C.~Dorlas,  A.~Ekert,  A.~Kay,  and A.~J.~Landahl, \emph{Perfect transfer of arbitrary states in quantum spin networks}, Phys.\ Rev.\ A, \textbf{71} (2005), 032312.
\bibitem{signed} J.~Brown, C.~Godsil, D.~Mallory, A.~Raz, and C.~Tamon, \emph{Perfect state transfer on signed graphs}, Quantum Inf.\ Comput.\ \textbf{13}, No.\ 5\& 6 (2013), 0511--0530.
\bibitem{Steve2015} S.\ Kirkland, \emph{Sensitivity analysis of perfect state transfer in quantum spin networks}, Linear Algebra  Appl.\ \textbf{472} (2015), 1--30.
\bibitem{BA} N.~Behzadi  and B.~Ahansaz, \emph{Enhancing quantum state transfer efficiency in binary-tree spin networks by partially collapsing measurements},  arXiv preprint \emph{arXiv:1611.03035}.
\bibitem{CRC} P.~Cappellaro, C.~Ramanathan,   and D.G.~Cory,  \emph{Simulations of information transport in spin chains}, Phys.\ Rev.\ Lett.\ \textbf{99}(25) (2007), 250506.
\bibitem{Kay15} A.~Kay, \emph{Quantum Error Correction for Noisy Quantum Wires}, arviv:1507.06139 (2015)


\bibitem{Pan} Y.~Pan, Z.~Miao, N.H.~Amini, V.~Ugrinovskii, and M.R.~James. \emph{Interpolation approach to Hamiltonian-varying quantum systems and the adiabatic theorem.} EPJ Quantum Technology \textbf{2}(1) (2015), 24.
\bibitem{Hamiltonian} M.~Christandl, N.~Datta, A.~Ekert,  and A.J.~Landahl, \emph{Perfect state transfer in quantum spin networks}. Phys.\ Rev.\ Lett., \textbf{92}(18) (2004), 187902.





\bibitem{BGS08} A.~Bernasconi, C.~Godsil, and S.~Severini, \emph{Quantum networks on cubelike graphs},
Phys.\ Rev.\ A \textbf{78} (2008), 052320.
\bibitem{C.Galgg} C.~Godsil and G.~Royle, \emph{Algebraic Graph Theory}. Springer, 2001.
\bibitem{GM} C.D.~Godsil and B.D.~McKay, \emph{Constructing cospectral graphs}, Aequat.\ Math.\ \textbf{25} (1982),
257-268.
\bibitem{ABWH} A.E.~Brouwer and W.H.~Haemers, \emph{Spectra of Graphs}. Springer, New York, 2012.


%\bibitem{quotientgraphpst} R.~Bachman, E.~Fredette, J.~Fuller, M.~Landry, M.~Opperman, C.~Tamon, and A.~Tollefson, \emph{Perfect state transfer of quantum walks on quotient graphs}, Quantum Information and Computation, \textbf{12}(3\& 4) (2012), 293-313.



%\bibitem{Sequencebook} C.\ Ding, T.\ Helleseth, and H.\ Niederreiter, \emph{Sequences and their Applications}, Springer-verlag, London, 1999.
\bibitem{Fiedleral} M.~Fiedler, \emph{Algebraic connectivity of graphs}, Czech.\ Math.\ J., \textbf{23}(98) (1973), 298-305.
\bibitem{quogra} R.\ Bachman, E.\ Fredette, J.\ Fuller, M.\ Landry, M.\ Opperman, C.\ Tamon, A.\ Tollefson, \emph{Perfect state transfer on quotient graphs}, Quantum Inf.\ Comput.\ \textbf{12} (2012)  293-313.
\bibitem{whenPST} C.~Godsil, \emph{When can perfect state transfer occur?}, Electronic Journal of Linear Algebra, \textbf{23}, (2012)


\bibitem{PRK} P.J.~Pemberton-Ross  and A.~Kay, \emph{Perfect quantum routing in regular spin networks}, Phys.\ Rev.\ Lett.\ \textbf{106}(2) (2011), 020503.

\bibitem{routingNik} G.M.~Nikolopoulos, \emph{Directional coupling for quantum computing and communication}, Phys.\ Rev.\ Lett.\ \textbf{101}(20) (2008), 200502.

\bibitem{routingCS} C.~Chudzicki and F.W.~Strauch, \emph{Parallel state transfer and efficient quantum routing on quantum networks}, Phys.\ Rev.\ Lett.\  \textbf{105}(26) (2010), 260501.

%\bibitem{HJ} R.A.~Horn and C.R.~Johnson, \emph{Topics in Matrix Analysis}. Cambridge university press, New York, 1991.

\bibitem{Dan} D.~Liberzon. \emph{Switching in systems and control}. Springer Science \& Business Media, 2012.

%\bibitem{BM} B.\ Mohar, \emph{The Laplacian spectrum of graphs}, Graph Theory, Combinatorics, and Applications, Vol.\ 2, Ed.\ Y.\ Alavi, G.~Chartrand, O.R.\ Oellermann, A.J.\ Schwenk, Wiley, 1991,  871-898.

 
\end{thebibliography}
\end{document}